\documentclass[12pt]{article}

\usepackage{amsmath,amssymb,amsfonts,amsthm,epic}

\usepackage{graphicx}

\newcommand\torder{\vartriangleright}

\theoremstyle{plain}
\newtheorem{thm}[subsubsection]{Theorem}

\newtheorem{ex}[subsubsection]{Example}

\theoremstyle{definition}
\newtheorem{rem}[subsubsection]{Remark}
\newtheorem{defn}[subsubsection]{Definition}
\newtheorem{nota}[subsubsection]{Notation}

\newcommand\brem{\begin{rem}\begin{sffamily}\begin{upshape}}
\newcommand\erem{\end{upshape}\end{sffamily}\end{rem}}
\newcommand\bdefn{\begin{defn}\begin{rm}}
\newcommand\edefn{\end{rm}\hfill$\Box$\end{defn}}
\newcommand\bnot{\begin{nota}\begin{rm}}
\newcommand\enot{\end{rm}\hfill$\Box$\end{nota}}
\newcommand\bex{\begin{ex}\begin{rm}}
\newcommand\eex{\end{rm}\hfill$\Box$\end{ex}}

\newenvironment{Proof}{%
\par\noindent{\scshape Proof:}\begin{rm}}{\hfill$\Box$\end{rm}\newline}
\numberwithin{equation}{subsection}
\title{An algorithm for the HK function for disjoint-term trinomial hypersurfaces}
\date{}
\author {Shyamashree Upadhyay\\Department of Mathematics\\ Indian Institute of Technology, Guwahati\\Guwahati, Assam-781039, INDIA\\shyamashree@iitg.ernet.in\\shyamashree.upadhyay@gmail.com}

\begin{document}
\maketitle

\begin{abstract}
A `trinomial hypersurface' is a hypersurface that is defined by a single polynomial having $3$ non-constant terms in it and no constant term. A `disjoint-term trinomial hypersurface' is a trinomial hypersurface whose defining polynomial has the property that any $2$ distinct terms in it have GCD equal to $1$. In this article, I provide an algorithm for computing the Hilbert-Kunz function for any disjoint-term trinomial hypersurface in general, over any field of arbitrary positive characteristic. However, I do not provide any formula for the Hilbert-Kunz function.
\end{abstract}
\tableofcontents
\section{Introduction}\label{s.Introduction}
Let $(A,\mathfrak{n})$ be a noetherian local ring of dimension $d$ and of prime characteristic $p>0$. Let $I$ be an $\mathfrak{n}$-primary ideal. Then the `Hilbert-Kunz function' of $A$ with respect to $I$ is defined as 
$$HK_{I,A}(p^n)=l(A/I^{(p^n)})$$ 
where $I^{(p^n)}$ = $n$-th Frobenius power of $I:=$ the ideal generated by $p^n$-th powers of elements of $I$.

The associated Hilbert-Kunz multiplicity is defined to be $$c(I,A)=\lim_{n\rightarrow\infty}\frac{HK_{I,A}(p^n)}{p^{nd}}$$.

Let $q$ denote an arbitrary positive power of $p$. Paul Monsky had proved in his paper \cite{Mo} that $$HK_{I,A}(q)=c(I,A)q^d+O(q^{d-1})$$ where $c(I,A)$ is a real constant. 

In many cases, it has been proved that the Hilbert-Kunz multiplicity is rational, see for example \cite{Conca}, \cite{Buch},\cite{Han}, \cite{Mo4}, \cite{Su} and \cite{Trivedi}. 

However, Paul Monsky has suggested in his paper \cite{Mo1} that modulo a conjecture, a certain hypersurface defined by a $5$-variable polynomial has irrational Hilbert-Kunz multiplicity. The total number of non-constant terms in this $5$-variable polynomial considered by Paul Monsky is $4$. Paul Monsky has few more papers in the same line (see for example \cite{Mo2} and \cite{Mo3}).

In this article, I am working on some special kind of hypersurfaces called `disjoint-term trinomial hypersurfaces'. I believe  that for showing that the Hilbert-Kunz multiplicity can become irrational in some cases, it is enough to work with such hypersurfaces, those defined by polynomials having $4$ non-constant terms in it (as in the example taken by Paul Monsky in his paper \cite{Mo1}) are not needed.

A `disjoint-term trinomial hypersurface' is defined in definition \ref{defn.disj-term-trinomialHyp} below. In the present work, I give an algorithm for computing the Hilbert-Kunz function for any disjoint-term trinomial hypersurface in general, over any field of arbitrary positive characteristic. And from this algorithm, I make a prediction that there can be examples of disjoint-term trinomial hypersurfaces for which the corresponding Hilbert-Kunz multiplicity can be irrational. I provide supportive reasoning for this prediction. More work supportive of this prediction is in progress. 

Before thinking about trinomial hypersurfaces, I have worked on the Hilbert-Kunz function for Binomial hypersurfaces (see \cite{Su}). A `binomial hypersurface' is a hypersurface that is defined by a single polynomial having $2$ non-constant terms in it and no constant term. In \cite{Su}, whatever is done before section $4$, holds true for any hypersurface over any field of positive characteristic, need not have to be a binomial hypersurface! A process called `mutation' defined in section $3$ of \cite{Su} which can be applied to any hypersurface in general. In fact, the process `mutation' becomes more rich in its combinatorial nature if used for hypersurfaces which are defined by polynomials having more than $2$ terms in it.

In this article, my strategy of work is to first apply the process `mutation' as mentioned above on trinomial hypersurfaces to reduce the problem to a problem of solving certain systems of linear equations (for a fixed trinomial hypersurface, we will need to solve $p^M$-many system of linear equations where $p$ is a prime number and $M$ is some positive integer). Then by applying some basic linear algebra techniques on the resulting problem, we are reduced to a problem of computing the ranks of a huge collection of matrices (which are of size of the order of $p^N$ for some positive integer $N$). This will produce the algorithm for computing the required Hilbert-Kunz function, which appears in section \ref{s.the-algorithm-irrational} of this article. This article gives this algorithm only, not the Hilbert-Kunz function. More work is in progress regarding how to use the algorithm to get the required Hilbert-Kunz function.

\section{Stating the problem}\label{s.statetheproblem}
Let $S=k[x_1,\cdots,x_m]$ where $k$ is a field of arbitrary prime characteristic $p>0$, and $J$ be an arbitrary ideal in $S$. Let $\mathfrak{m}$ be the maximal ideal $(x_1,\cdots,x_m)$ of $S$ and let $R=S/J$. Then $\hat{\mathfrak{m}}:=\mathfrak{m}+J$ is a maximal ideal in $R$. Without loss of generality, we can assume that $J\subseteq\mathfrak{m}$, for otherwise $\hat{\mathfrak{m}}$ is the whole ring $R$. The `Hilbert-Kunz function' of the noetherian local ring $R_{\hat{\mathfrak{m}}}$ with respect to $\hat{\mathfrak{m}}R_{\hat{\mathfrak{m}}}$ is given by:
$$HK_{\hat{\mathfrak{m}}R_{\hat{\mathfrak{m}}},R_{\hat{\mathfrak{m}}}}(p^n)=l(\frac{R_{\hat{\mathfrak{m}}}}{(\hat{\mathfrak{m}}R_{\hat{\mathfrak{m}}})^{(p^n)}})$$
where $(\hat{\mathfrak{m}}R_{\hat{\mathfrak{m}}})^{(p^n)}$ = $n$-th Frobenius power of $\hat{\mathfrak{m}}R_{\hat{\mathfrak{m}}}$.

Note that the rings $\frac{R_{\hat{\mathfrak{m}}}}{(\hat{\mathfrak{m}}R_{\hat{\mathfrak{m}}})^{(p^n)}}$ and $\frac{R}{{\mathfrak{m}}^{(p^n)}+J}$ are isomorphic, where ${\mathfrak{m}}^{(p^n)}=(x_1^{p^n},\cdots,x_m^{p^n})$. So it is enough to compute the length $l(\frac{R}{{\mathfrak{m}}^{(p^n)}+J})$, i.e., 
$$HK_{\hat{\mathfrak{m}}R_{\hat{\mathfrak{m}}},R_{\hat{\mathfrak{m}}}}(p^n)=l(\frac{R}{{\mathfrak{m}}^{(p^n)}+J})$$
This function is called the \textit{Hilbert-Kunz function of $R$} (with respect to $\hat{\mathfrak{m}}$).  

We are interested in computing the Hilbert-Kunz function of $R$ when $R$ is a `disjoint-term trinomial hypersurface'. A `disjoint-term trinomial hypersurface' is defined as follows:
\bdefn\label{defn.disj-term-trinomialHyp} 
Let $S$ and $J$ be as above. Let $J=(f)$ where $f$ is a polynomial that contains $3$ non-constant terms and no constant term. The affine variety defined by the ideal $J$ is called a \textbf{Trinomial Hypersurface}. A `disjoint-term trinomial hypersurface' is a trinomial hypersurface whose defining polynomial has the property that any $2$ distinct terms in it have GCD equal to $1$.
\edefn

\subsection{The filtration for computing the length}\label{ss.filtration-and-checking}
Recall from \S 2.1 of \cite{Su} the filtration for computing the length $l(\frac{R}{{\mathfrak{m}}^{(p^n)}+J})$ which was having the property that \textit{each successive quotient} is \textit{either} \textbf{one dimensional} \textit{or} \textbf{zero}. From \S 2.2 of \cite{Su}, we can also recall the \textbf{key checking} which needs to be done at each step for computing the length $l(\frac{R}{{\mathfrak{m}}^{(p^n)}+J})$. This checking was the following:
\begin{equation}\label{eq.key-check}
a_t\in I_{t-1}+J\ \text{or not for every}\ t\in\{1,\ldots,p^{mn}\}
\end{equation}
where $a_t$ and $I_{t-1}$ are as defined in \S 2.1 of \cite{Su}. It is easy to see that for any $t\in\{1,\ldots,p^{mn}\}$, if $a_t\in I_{t-1}+J$, then the quotient $\frac{I_t+J}{I_{t-1}+J}$ is \textit{zero}, and the quotient $\frac{I_t+J}{I_{t-1}+J}$ is \textit{one-dimensional} otherwise. The length $l(\frac{R}{{\mathfrak{m}}^{(p^n)}+J})$ is the cardinality of the set $\{t|t\in\{1,\ldots,p^{mn}\}$ for which the quotient $\frac{I_t+J}{I_{t-1}+J}$ is one-dimensional$\}$. 
\section{Reduction to linear algebra}\label{s.red-lin-alg}
In this section, we will first define a term order $\torder$ on the set of all monomials in the variables $x_1,\ldots,x_m$, and then with respect to $\torder$, we will arrange the terms of the polynomial $f$. With the help of all this notation, we will apply the process `mutation' (which was introduced in \S 3.2 of \cite{Su}) to do the checking \ref{eq.key-check}. Using this process, we will then transform the original problem to a linear algebra problem.
\subsection{The term order $\torder$}\label{ss.torder}
Let us put an order (denote it by $\torder$) on the set of all monomials in the variables $x_1,\ldots,x_m$ as follows:---
\begin{itemize}
\item Set $x_1\torder\cdots\torder x_m$.
\item On the set of all monomials in the variables $x_1,\ldots,x_m$, $\torder$ is the degree lexicographic order with respect to the order $\torder$ defined on the variables $x_1,\ldots,x_m$.
\end{itemize} 
Since the polynomial $f$ has $3$ non-constant terms in it, let us denote the most initial (with respect to $\torder$) term of $f$ as $[3]$, the next most initial term of $f$ as $[2]$ and the least initial term as $[1]$. Hence we have
$$J=(f)=([3]+[2]+[1])$$
\brem\label{r.scalarsintermsoff}
Note here that the terms $[3],[2],[1]$ of $f$ are assumed to be containing scalar coefficients.
\erem  
\subsection{The process Mutation}\label{ss.mutation}
Recall the set $\mathfrak{M}:=\{a_t|t\in\{1,\ldots,p^{mn}\}\}$ from \S 2.1 of \cite{Su}. Let $A$ be an arbitrary element of the set $\mathfrak{M}$. Say the monomial $A$ equals $a_l$ for some $l\in\{1,\ldots,p^{mn}\}$. The key-checking condition for the monomial $A$ says that `$A\in I_{l-1}+J$ or not'. Let us denote by $A_c$ the ideal $I_{l-1}$. We call $A_c$ the \textbf{ key ideal corresponding to A}.
Let us call any monomial (not equal to $A$) which belongs to the ideal $A_c$ as a \textit{convergent term with respect to the ideal $A_c$}, and any monomial which does not belong to the ideal $A_c$ as a \textit{non-convergent term with respect to the ideal $A_c$}. For convenience of terminology, we will henceforth omit the phrase `with respect to the ideal $A_c$', unless otherwise needed, and continue calling a monomial to be \textit{convergent} or \textit{non-convergent}. 

Recall the process `mutation' and all remark(s), definitions and notation related to it from \S 3.2 of \cite{Su}. Theorem 3.2.8 of \cite{Su} can be restated for the present case of disjoint-term trinomial hypersurfaces in the following way:
\begin{thm}\label{t.mutation-stop}
Let $f=[3]+[2]+[1]$. $A\in A_c+J$ if and only if one of the following mutually exclusive conditions hold:---\\
(i) The term $[1]$ of $f$ divides the monomial $A$.\\
\noindent (ii) The term $[1]$ of $f$ does not divide the monomial $A$, but there exists term(s) of $f$ not equal to $[1]$ which divide $A$ and the mutation process (with respect to the monomial $A$ and the polynomial $f$) stops. 
\end{thm}
Given any term $\tau$ of the polynomial $f$, define $[-\tau]:=\frac{1}{\tau}$. Let
\begin{center}
$\mathcal{B}(A,f):=$ the set of all non-convergent mutants in $A$ and $f$ and\\
$\mathcal{A}(A,f):=\{B[-\tau]|B\in\mathcal{B}(A,f),\tau\neq 1\ and\ \tau\ divides\ B\}$. 
\end{center}
Clearly every element of $\mathcal{A}(A,f)$ is a mutator in $A$ and $f$. With all this notation, we now modify theorem \ref{t.mutation-stop} above so that the resulting theorem becomes more convenient for our purpose of doing the checking as mentioned in equation \ref{eq.key-check}. The following theorem provides the modification:
\begin{thm}\label{t.mutation-stop-modified}
Let $f=[3]+[2]+[1]$. $A\in A_c+J$ if and only if one of the following mutually exclusive conditions hold:---\\
(i) The term $[1]$ of $f$ divides the monomial $A$.\\
\noindent (ii) The term $[1]$ of $f$ does not divide the monomial $A$, but there exist term(s) of $f$ not equal to $[1]$ which divide $A$ and there exists scalars $c_D$ (corresponding to each $D\in\mathcal{A}(A,f)$) such that the product $f.(\Sigma_{D\in\mathcal{A}(A,f)}c_{D}D)$ equals $bA+finitely\ many\ convergent\ terms$ where $b$ is a non-zero scalar.
\end{thm}
\begin{Proof}
The proof is immediate from theorem \ref{t.mutation-stop} and the definition of the process mutation.
\end{Proof}
We can now modify theorem \ref{t.mutation-stop-modified} by separating out the case when $[2]$ divides the monomial $A$ and there exists a positive integer $M$ for which $A[-2]^M$ contains no negative powers and $A\dfrac{{[-2]}^M}{{[1]}^{M}}$ is convergent (note here that by the notation $[-2]^M$, we mean $[-2]$ multiplied $M$ times, and similarly for the notation $[1]^M$). The modified version is the following:
\begin{thm}\label{t.mutation-stop-modifiedv2}
Let $f=[3]+[2]+[1]$. $A\in A_c+J$ if and only if one of the following mutually exclusive conditions hold:---\\
(i) The term $[1]$ of $f$ divides the monomial $A$.\\
\noindent (ii) The term $[1]$ of $f$ does not divide the monomial $A$, but the term $[2]$ divides $A$ and there exists a positive integer $M$ for which $A[-2]^M$ contains no negative powers and $A\dfrac{{[-2]}^M}{{[1]}^{M}}$ is convergent. \\
\noindent (iii) Neither condition (i) nor condition (ii) holds, but $[3]$ divides $A$ ($[2]$ may or may not divide $A$, but if at all $[2]$ divides $A$, condition (ii) does not hold) and there exists scalars $c_D$ (corresponding to each $D\in\mathcal{A}(A,f)$) such that the product $f.(\Sigma_{D\in\mathcal{A}(A,f)}c_{D}D)$ equals $bA+finitely\ many\ convergent\ terms$ for some non-zero scalar $b$.
\end{thm}
\begin{Proof}
It follows from the construction of the ideal $A_c$ and from the disjointness of the terms of the polynomial $f$ that any mutant in $A$ and $f$ which when expressed in lowest terms contains only $[-2]$s in the numerator and at least one $[3]$ in the denominator, is convergent. This fact together with theorem \ref{t.mutation-stop-modified} yield the result: Let $M_0$ be the smallest positive integer $M$ which satisfies condition (ii) of this theorem. Multiply $f$ by a suitable scalar linear combination of mutators of the form $A\dfrac{{[-2]}^n}{{[1]}^{n-1}}$ where $1\leq n\leq M_0$. 
\end{Proof}
We will now provide an equivalent formulation of condition (iii) of theorem \ref{t.mutation-stop-modifiedv2} above in terms of linear algebra. This equivalent formulation is condition $(iii)'$ of remark \ref{r.condition(iii)-to-sysoflineareqn} below.
\brem\label{r.condition(iii)-to-sysoflineareqn}
Let $\mathcal{E}(A,f):=\{D[\tau]|D\in\mathcal{A}(A,f)\ and\ \tau\in\{[1],[2],[3]\}\}$ and $\mathcal{L}(A,f):=$the set of all elements in $\mathcal{E}(A,f)$ which when expressed in lowest terms do not contain any $[3]$ in the denominator. Note that the product $f.(\Sigma_{D\in\mathcal{A}(A,f)}c_{D}D)$ in condition (iii) of theorem \ref{t.mutation-stop-modifiedv2} above equals a linear combination of elements of $\mathcal{E}(A,f)$, say $\Sigma_{B\in\mathcal{E}(A,f)}e_{B}B$. This sum can be broken into $2$ parts as follows: $\Sigma_{B\in\mathcal{E}(A,f)}e_{B}B=\Sigma_{B\notin\mathcal{L}(A,f)}e_{B}B+\Sigma_{B\in\mathcal{L}(A,f)}e_{B}B$. Since any mutant in $A$ and $f$ which when expressed in lowest terms contains only $[-2]$s in the numerator and at least one $[3]$ in the denominator is convergent, it follows that the portion $\Sigma_{B\notin\mathcal{L}(A,f)}e_{B}B$ contains all convergent terms and they are finitely many. So if we equate the coefficients of like terms of the product $f.(\Sigma_{D\in\mathcal{A}(A,f)}c_{D}D)$ and the sum $\Sigma_{B\in\mathcal{L}(A,f)}e_{B}B$, we get a system $\mathfrak{A}_{A,f}\mathfrak{X}=\mathfrak{B}$ of linear equations where $\mathfrak{X}$ is a column vector in the unknowns $c_D$, $\mathfrak{B}$ is a column vector in the scalars $e_B$ where $B\in\mathcal{L}(A,f)$ and $\mathfrak{A}_{A,f}$ is a matrix with entries from the set $\{0,1\}$. 
\erem
So an equivalent formulation of condition (iii) of theorem \ref{t.mutation-stop-modifiedv2} above will be condition  $(iii)'$ as stated  below: 
\begin{quote}
$(iii)'$ Neither condition (i) nor condition (ii) of theorem \ref{t.mutation-stop-modifiedv2} holds, but $[3]$ divides $A$ and the system $\mathfrak{A}_{A,f}\mathfrak{X}=\mathfrak{B}$ of linear equations is solvable for the vector $\mathfrak{B}$ which is having the property that:
\begin{center}
$e_B\neq 0$ for $B=A$ and $e_B=0$ for all $B$ non-convergent.
\end{center}
\end{quote}
Combining theorem \ref{t.mutation-stop-modifiedv2} and condition $(iii)'$ above, we get the following theorem:
\begin{thm}\label{t.mutation-stop-modifiedv3}
$A\in A_c+J$ if and only if 
\begin{itemize}
\item either condition (i) or condition (ii) of theorem \ref{t.mutation-stop-modifiedv2} holds or
\item neither condition (i) nor (ii) of theorem \ref{t.mutation-stop-modifiedv2} holds but $[3]$ divides $A$ and the system $\mathfrak{A}_{A,f}\mathfrak{X}=\mathfrak{B}$ of linear equations is solvable for the vector $\mathfrak{B}$ (this vector is introduced in remark \ref{r.condition(iii)-to-sysoflineareqn} above) which is having the property that:
\end{itemize}
\begin{center}
$e_B\neq 0$ for $B=A$ and $e_B=0$ for all $B$ non-convergent.
\end{center}
\end{thm}
\subsection{Further reduction to another linear system}\label{ss.further-red-toanother}
To simplify the problem, we will now transform the system $\mathfrak{A}_{A,f}\mathfrak{X}=\mathfrak{B}$ of linear equations further to an equivalent system $\mathfrak{A}_{A,f}'\mathfrak{X}=\mathfrak{B}'$ of linear equations. But for that we first need to introduce some definitions and notation.
\bdefn\label{defn.num-den-ofmutator}
Let $D\in\mathcal{A}(A,f)$. After expressing $D$ in its lowest terms, consider the product of all terms of the form $[-\tau]$ (where $\tau\in\{[2],[3]\}$) appearing in the lowest-term expression of $D$. This product is called the \textit{numerator} of $D$. Similarly the product of all terms of the form $[\tau]$ (where $\tau\in\{[1],[2],[3]\}$) appearing in the lowest-term expression of $D$ is called the \textit{denominator} of $D$.\\
A similar idea applies to the definition(s) of \textit{numerator} and \textit{denominator} of any element of $\mathcal{L}(A,f)$. 
\edefn
\bdefn\label{defn.standard-form-mutator}
Let $D\in\mathcal{A}(A,f)$. Express $D$ in its lowest terms and consider the numbers appearing in the numerator of $D$ ignoring the $-$ sign. Arrange these numbers to get the lowest possible positive integer. We call this integer as the \textit{standard form of the numerator} of $D$ and denote it by $D_{num}$. \\
Similarly after expressing $D$ in its lowest terms, one can consider the numbers appearing in the denominator of $D$ (of course, no integer comes with a $-$ sign in the denominator). Arrange these numbers to get the lowest possible positive integer. We call this integer as the \textit{standard form of the denominator} of $D$ and denote it by $D_{den}$. 
\edefn
\bnot\label{n.preceq}
Let us put an order $\preceq$ on elements of $\mathcal{A}(A,f)$. For any two elements $D,D'\in\mathcal{A}(A,f)$, we say that $D\preceq D'$ if exactly one of the following $2$ conditions hold:
\begin{itemize}
 \item $D_{num}<D'_{num}$.
 \item $D_{num}=D'_{num}$ but $D_{den}\leq D'_{den}$. 
\end{itemize}
Clearly, $\preceq$ is a total order on the elements of $\mathcal{A}(A,f)$.
\enot
\bdefn\label{defn.standard-form-mutant}
Let $B\in\mathcal{L}(A,f)$. Express $B$ in its lowest terms and consider the numbers appearing in the numerator of $B$ ignoring the $-$ sign. Arrange these numbers to get the lowest possible positive integer. We call this integer as the \textit{standard form of the numerator} of $B$ and denote it by $B_{num}$. \\
Similarly after expressing $B$ in its lowest terms, one can consider the numbers appearing in the denominator of $B$ (of course, no integer comes with a $-$ sign in the denominator). Arrange these numbers to get the lowest possible positive integer. We call this integer as the \textit{standard form of the denominator} of $B$ and denote it by $B_{den}$. 
\edefn
\bnot\label{n.precsim}
Let us put an order $\precsim$ on elements of $\mathcal{L}(A,f)$. For any two elements $B,B'\in\mathcal{L}(A,f)$, we say that $B\preceq B'$ if exactly one of the following $2$ conditions hold:
\begin{itemize}
 \item $B_{num}<B'_{num}$.
 \item $B_{num}=B'_{num}$ but $B_{den}\leq B'_{den}$. 
\end{itemize}
Clearly, $\precsim$ is a total order on the elements of $\mathcal{L}(A,f)$.
\enot
Let $A$ be a monomial such that neither condition (i) nor condition (ii) of theorem \ref{t.mutation-stop-modifiedv2} hold. Consider the system $\mathfrak{A}_{A,f}\mathfrak{X}=\mathfrak{B}$ of linear equations as mentioned in condition $(iii)'$ above. Recall that $\mathfrak{X}$ is a column vector in the unknowns $c_D$ where $D\in\mathcal{A}(A,f)$, $\mathfrak{B}$ is a column vector in the scalars $e_B$ where $B\in\mathcal{L}(A,f)$. Let us arrange the unknowns $c_D$ of the vector $\mathfrak{X}$ in the `increasing' order induced by the order $\preceq$ on $\mathcal{A}(A,f)$. But we will arrange the elements of the vector $\mathfrak{B}$ differently and for describing this arrangement, we need the following notation:
\bnot\label{n.partition-of-L(A,f)}
Let $1_{min,A}$ (resp. $2_{min,A}$) be the minimum number of $[1]$s (resp. $[2]$s) needed to be multiplied to the monomial $A$ for convergence. Let $(-2)_{max,A}$ (resp. $(-3)_{max,A}$) be the maximum number of $2$s (resp. $3$s) that can be divided from $A$ (in such a way that the quotient does not contain negative power of any of the underlying variables). 
\begin{center}
\underline{Case I: When $(-2)_{max,A}\geq 1$}
\end{center} 
Let $\mathcal{R}(A,f):=\mathcal{L}(A,f)\cap(\{A\dfrac{[-2]}{[1]}\}\cup\{A\dfrac{[-3]^i}{[1]^{i-1}[2]}|i=1,\ldots,1_{min,A}+1\}\cup\{A\dfrac{[-3]^j}{[1]^{1_{min,A}}[2]^k}|2\leq k\leq 2_{min,A}-1\ and\ j=1_{min,A}+k\})$. 
\begin{center}
\underline{Case II: When $(-2)_{max,A}=0$} 
\end{center}
Let $\mathcal{R}(A,f):=\mathcal{L}(A,f)\cap(\{A\dfrac{[-3]}{[1]}\}\cup\{A\dfrac{[-3]^i}{[1]^{i-1}[2]}|i=1,\ldots,1_{min,A}+1\}\cup\{A\dfrac{[-3]^j}{[1]^{1_{min,A}}[2]^k}|2\leq k\leq 2_{min,A}-1\ and\ j=1_{min,A}+k\})$. 
\begin{center}
\underline{For cases I and II both}
\end{center} 
Let $\mathcal{S}(A,f):=\mathcal{L}(A,f)-\mathcal{R}(A,f)$.
\enot
\noindent The arrangement for the vector $\mathfrak{B}$:
\begin{quote}
We will first consider those elements of $\mathfrak{B}$ which correspond to the elements of $\mathcal{S}(A,f)$ and arrange them in the `increasing' order induced by the order $\precsim$ on $\mathcal{L}(A,f)$. Below that, we put those elements of $\mathfrak{B}$ which correspond to the elements of $\mathcal{R}(A,f)$ and arrange them in the `decreasing' order induced by the order $\precsim$ on $\mathcal{L}(A,f)$. 
\end{quote}
\noindent With this arrangement of the elements of the vectors $\mathfrak{X}$ and $\mathfrak{B}$, the matrix $\mathfrak{A}_{A,f}$ achieves the following property in its structure:
\begin{quote}
The submatrix of $\mathfrak{A}_{A,f}$ corresponding to the first $|\mathcal{S}(A,f)|$-many rows of it is an upper triangular matrix having $1$-s on the diagonal and all entries above the diagonal lie in the set $\{0,1\}$. Let $T(A,f)$ denote this upper triangular matrix.
\end{quote}
Now using some row operations, we can transform the system $\mathfrak{A}_{A,f}\mathfrak{X}=\mathfrak{B}$ of linear equations (arranged in the above mentioned fashion) into another equivalent system $\mathfrak{A}_{A,f}'\mathfrak{X}=\mathfrak{B}'$ of linear equations where the matrix $\mathfrak{A}_{A,f}'$ and the vector $\mathfrak{B}'$ are given by the following description:\\
\begin{quote}
The submatrix of $\mathfrak{A}_{A,f}'$ corresponding to the first $|\mathcal{S}(A,f)|$-many rows of it is the same as the submatrix of $\mathfrak{A}_{A,f}$ corresponding to the first $|\mathcal{S}(A,f)|$-many rows of it and, the rest of the matrix $\mathfrak{A}_{A,f}'$ is the zero matrix.

The matrix $\mathfrak{A}_{A,f}'$ is obtained from the matrix $\mathfrak{A}_{A,f}$ by applying some row operations: One can make any one of the last $|\mathcal{R}(A,f)|$-many rows of $\mathfrak{A}_{A,f}$ zero by using the upper triangular submatrix $T(A,f)$ of $\mathfrak{A}_{A,f}$ and applying appropriate row operations. For instance, say $R$ is any row among the last $|\mathcal{R}(A,f)|$-many rows of $\mathfrak{A}_{A,f}$. Consider the first (from left to right) entry in $R$ which is non-zero, say this happens at column $C$ and let $c$ be this non-zero entry. Go to that row of the upper triangular matrix $T(A,f)$ whose first (from left to right) non-zero entry is at column $C$, say this happens at row $R'$. Replace row $R$ by $R-cR'$. In the resulting row, look at the first non-zero entry. Then proceed similarly using the matrix $T(A,f)$ and appropriate row operations till row $R$ (that we started with) becomes zero. 

The column vector $\mathfrak{B}'$ is obtained from the vector $\mathfrak{B}$ by applying the same row operations which were used to obtain the matrix $\mathfrak{A}_{A,f}'$ from the matrix $\mathfrak{A}_{A,f}$. 
\end{quote}
\section{Reduction to Combinatorics}\label{s.red-to-comb}
In this section, we will reduce the problem further to a problem of solving another equivalent system $\mathfrak{B}_{A,f}\mathfrak{Y}=\mathfrak{e}$ of linear equations, which has nice combinatorial properties.
\subsection{Description of a new system}\label{ss.description-of-anewsys}
The system $\mathfrak{A}_{A,f}'\mathfrak{X}=\mathfrak{B}'$ of linear equations as mentioned above is solvable if and only if the entries in the last $|\mathcal{R}(A,f)|$-many rows of the vector $\mathfrak{B}'$ are all zero. Equating the entries in the last $|\mathcal{R}(A,f)|$-many rows of the vector $\mathfrak{B}'$ to zero, we get a new system $\mathfrak{B}_{A,f}\mathfrak{Y}=\mathfrak{e}$ of linear equations which can be described as follows:\\
\begin{itemize}
\item The number of rows in the matrix $\mathfrak{B}_{A,f}$ is $|\mathcal{R}(A,f)|$. The rows of the matrix $\mathfrak{B}_{A,f}$ are indexed by elements of the set $\mathcal{R}(A,f)$ which are arranged in the `decreasing' order induced by the order $\precsim$ on $\mathcal{L}(A,f)$.
\item The number of columns in it is $|P_{A,f}|$ where $P_{A,f}:=$ the set of all convergent elements in $\mathcal{L}(A,f)$. The columns of the matrix $\mathfrak{B}_{A,f}$ are indexed by elements of the set $P_{A,f}$.
\item In the case when $(-2)_{max,A}\geq 1$, the column vector $\mathfrak{e}$ is given by $[0,\ldots,0,e_A]^t$ where $e_A\neq 0$. And in the case when $(-2)_{max,A}=0$, the column vector $\mathfrak{e}$ is given by $[0,\ldots,0,e_A,e_A]^t$ where $e_A\neq 0$.
\item The entries of the column vector $\mathfrak{Y}$ are the elements of the vector $\mathfrak{B}$ which correspond to the elements of the set $P_{A,f}$ and they are arranged in the `increasing' order induced by the order $\precsim$ on $\mathcal{L}(A,f)$. 
\item The entries of the matrix $\mathfrak{B}_{A,f}$ have a nice combinatorial pattern which is mentioned in the two cases considered below. 
\end{itemize}
\begin{center}
\underline{Case I: When $(-2)_{max,A}\geq 1$}
\end{center}
\noindent\textbf{The last row of} $\mathfrak{B}_{A,f}$: Consider those elements of $P_{A,f}$ which when written in its lowest terms contains at least one $[2]$ in the denominator. The entry in the last row of $\mathfrak{B}_{A,f}$ corresponding to any such column is $0$. Next consider those elements of $P_{A,f}$ which when written in lowest terms do not contain any $[2]$ in the denominator. Such elements are of the form $A\dfrac{[-2]^i[-3]^j}{[1]^{1_{min,A}}}$ where $0\leq i,j\leq 1_{min,A}$, $i+j=1_{min,A}$ and $A\dfrac{[-2]^i[-3]^j}{[1]^{1_{min,A}}}$ contains no negative powers of any of the underlying variables. The entry in the last row of $\mathfrak{B}_{A,f}$ at the column corresponding to $A\dfrac{[-2]^i[-3]^j}{[1]^{1_{min,A}}}$ equals $(-1)^{(1_{min,A}+1)}\ ^{(1_{min,A})}C_i$. 

\noindent\textbf{Any row of} $\mathfrak{B}_{A,f}$ \textbf{other than the last row:} The columns of the matrix $\mathfrak{B}_{A,f}$ are indexed by elements of the set $P_{A,f}$. The elements of $P_{A,f}$ when written in their lowest terms exhibit the following pattern :
\begin{itemize}
\item \textbf{Type A:} Either it contains exactly $1_{min,A}$-many $[1]$s in the denominator and the number of $[2]$s in the denominator is strictly less than $2_{min,A}$. 
\item \textbf{Type B:} Or it contains exactly $2_{min,A}$-many $[2]$s in the denominator and the number of $[1]$s in the denominator is strictly less than $1_{min,A}$. 
\end{itemize}

First consider the \textit{type A} elements. These elements are either of the form $A\dfrac{[-2]^i[-3]^j}{[1]^{1_{min,A}}}$ where $0\leq i,j\leq 1_{min,A}$ and $i+j=1_{min,A}$ or $A\dfrac{[-3]^m}{[1]^{1_{min,A}}[2]^j}$ where $1\leq j\leq 2_{min,A}-1$ and $m=1_{min,A}+j$ where neither $A\dfrac{[-2]^i[-3]^j}{[1]^{1_{min,A}}}$ nor $A\dfrac{[-3]^m}{[1]^{1_{min,A}}[2]^j}$ contains any negative power of any of the underlying variables. Look at the column corresponding to the element $A\dfrac{[-2]^i[-3]^j}{[1]^{1_{min,A}}}$ of $P_{A,f}$. For any $0\leq r\leq 1_{min,A}-1$, consider the row corresponding to the entry $A\dfrac{[-3]^k}{[1]^{r}[2]}$ of $\mathcal{R}(A,f)$ where $k=r+1$. The entry in such a row of $\mathfrak{B}_{A,f}$ at this column equals $(-1)^{(1_{min,A}-k+1)}\ ^{(1_{min,A}-k)}C_i$. For any other row, the entry of $\mathfrak{B}_{A,f}$ at this column equals $0$. Next look at the column corresponding to the element $A\dfrac{[-3]^m}{[1]^{1_{min,A}}[2]^j}$ of $P_{A,f}$. There exists a row in $\mathfrak{B}_{A,f}$ which corresponds to the same element of $\mathcal{R}(A,f)$. The entry of $\mathfrak{B}_{A,f}$ in this column equals $1$ if the row corresponds to the same element and equals $0$ otherwise.\\

Next consider the \textit{type B} elements. These elements are of the form $A\dfrac{[-3]^m}{[1]^s[2]^{2_{min,A}}}$ where $0\leq s\leq 1_{min,A}-1$, $m$ is some integer such that $s+2_{min,A}=m$ and $A\dfrac{[-3]^m}{[1]^s[2]^{2_{min,A}}}$ contains no negative powers of any of the underlying variables. Look at the column corresponding to the element $A\dfrac{[-3]^m}{[1]^s[2]^{2_{min,A}}}$ of $P_{A,f}$. For any $0\leq r\leq 1_{min,A}-1$, consider the row corresponding to the entry $A\dfrac{[-3]^k}{[1]^{r}[2]}$ of $\mathcal{R}(A,f)$ where $k=r+1$. The entry in such a row of the matrix $\mathfrak{B}_{A,f}$ at this column equals $(-1)^{(2_{min,A+1})}\ ^{(m-1-s)}C_{r-s}$ if $s\leq r$ and $0$ if $s>r$. For any other row, it is not needed to know the entry of $\mathfrak{B}_{A,f}$ at the column corresponding to the element $A\dfrac{[-3]^m}{[1]^s[2]^{2_{min,A}}}$ of $P_{A,f}$, the reason behind this is explained in remark \ref{r.no-need-toknow-entries} below. 
\brem\label{r.no-need-toknow-entries}
We will see in section \ref{s.solvability-of-thenewsys} below that the system $\mathfrak{B}_{A,f}\mathfrak{Y}=\mathfrak{e}$ of linear equations can be transformed into an equivalent and simpler system $\mathfrak{B}_{A,f}^{eq}\mathfrak{Y}^{eq}=\mathfrak{e}^{eq}$ of linear equations where the matrix $\mathfrak{B}_{A,f}^{eq}$ is of the form $\begin{pmatrix}P & Q\\R & S\end{pmatrix}$ where $P,Q,R,S$ are blocks of the matrix $\mathfrak{B}_{A,f}^{eq}$ having the following properties:\\
\begin{center}
$P$ is a matrix having $1$ s or $-1$ s in the antidiagonal and $0$ s elsewhere, $R$ is the zero matrix, $Q,S$ are matrices of appropriate sizes with real entries.
\end{center}
It is elementary to see that for solving the system $\mathfrak{B}_{A,f}^{eq}\mathfrak{Y}^{eq}=\mathfrak{e}^{eq}$ of linear equations, it is not necessary to know the entries of the matrix $Q$. Let $l$ denote the total number of rows of $S$. Let $\mathfrak{Y}^{eq}_l$ and $\mathfrak{e}^{eq}_l$ denote the column vectors formed by the last $l$ entries of the column vectors $\mathfrak{Y}^{eq}$ and $\mathfrak{e}^{eq}$ respectively. In fact, solving the system $\mathfrak{B}_{A,f}^{eq}\mathfrak{Y}^{eq}=\mathfrak{e}^{eq}$ is equivalent to solving the system $S\mathfrak{Y}^{eq}_l=\mathfrak{e}^{eq}_l$. Therefore it is not needed to know certain entries of the matrix $\mathfrak{B}_{A,f}$.
\erem
\begin{center}
\underline{Case II: When $(-2)_{max,A}=0$} 
\end{center}
The matrix $\mathfrak{B}_{A,f}$ in this case differs from that in case I only at the following points:\\
\begin{itemize}
 \item There is no column of the matrix $\mathfrak{B}_{A,f}$ which is indexed by type A elements of the form $A\dfrac{[-2]^i[-3]^j}{[1]^{1_{min,A}}}$ where $1\leq i\leq 1_{min,A}$, $0\leq j\leq 1_{min,A}$ and $i+j=1_{min,A}$.
 \item In the last but one row (that is, the row corresponding to the entry $A\dfrac{[-3]}{[2]}$ of $\mathcal{R}(A,f)$), the entry at the column corresponding to the element $A[\dfrac{-3}{1}]^{1_{min,A}}$ is $0$. 
 \item The entry of the vector $\mathfrak{e}$ at the last but one-th row is $e_A$ which is $\neq 0$. 
\end{itemize}
Everything else remains the same as in case I. This finishes the description of the matrix $\mathfrak{B}_{A,f}$ when $(-2)_{max,A}=0$.
\brem\label{r.equivalent-matrix-incaseII}
When $(-2)_{max,A}=0$, the system $\mathfrak{B}_{A,f}\mathfrak{Y}=\mathfrak{e}$ of linear equations can be changed to an equivalent system by performing the following row operation:
\begin{center}
Replace the `last but one-th row' by `the last but one-th row \textit{minus} the last row'.
\end{center}
In the transformed system, the vector $\mathfrak{e}$ has also changed to the vector $[0,\ldots,0,e_A]^t$ where $e_A\neq 0$ (this is like the vector $\mathfrak{e}$ of case I). Let us denote the transformed system by $\mathfrak{B}_{A,f}^{tr}\mathfrak{Y}^{tr}=\mathfrak{e}^{tr}$. The index set of the columns of the matrix $\mathfrak{B}_{A,f}^{tr}$ remains the same as that of the matrix $\mathfrak{B}_{A,f}$ in the case when $(-2)_{max,A}\geq 1$. 
\erem
Due to a reason similar to that mentioned in remark \ref{r.no-need-toknow-entries} above (which was for case I), it is not needed to know some entries of the matrix $\mathfrak{B}_{A,f}^{tr}$.
\subsection{The theorem restated}\label{ss.the-thm-restated}
\brem\label{r.-2and-3maxgeq1min}
It follows easily from the description of the last row of the matrix $\mathfrak{B}_{A,f}$ that the system $\mathfrak{B}_{A,f}\mathfrak{Y}=\mathfrak{e}$ of linear equations is not solvable if $(-3)_{max,A}<1_{min,A}-(-2)_{max,A}$ or in other words, if $(-2)_{max,A}+(-3)_{max,A}<1_{min,A}$. The preceding statement is valid for the possibility $(-2)_{max,A}=0$ as well. 
\erem
\noindent We now have the following theorem:
\begin{thm}\label{t.mutation-stop-modifiedv4}
Let $f=[3]+[2]+[1]$. $A\in A_c+J$ if and only if one of the following mutually exclusive conditions hold:---\\
(i) The term $[1]$ of $f$ divides the monomial $A$.\\
\noindent (ii) The term $[1]$ of $f$ does not divide the monomial $A$, but the term $[2]$ divides $A$ and there exists a positive integer $M$ for which $A[-2]^M$ contains no negative powers and $A\dfrac{{[-2]}^M}{{[1]}^{M}}$ is convergent. \\
\noindent (iii) Neither condition (i) nor condition (ii) holds, but $[3]$ divides $A$, and\\
\begin{itemize}
\item Either $(-2)_{max,A}=0$, $(-3)_{max,A}\geq 1_{min,A}$ and the system $\mathfrak{B}_{A,f}^{tr}\mathfrak{Y}^{tr}=\mathfrak{e}^{tr}$ (as mentioned in remark ~\ref{r.equivalent-matrix-incaseII} above) of linear equations is solvable.
\item or $(-2)_{max,A}\geq 1$, $(-2)_{max,A}+(-3)_{max,A}\geq 1_{min,A}$ and the system $\mathfrak{B}_{A,f}\mathfrak{Y}=\mathfrak{e}$ of linear equations is solvable.
\end{itemize}
\end{thm}
\subsection{A simpler description of the new system}\label{ss.new-combinatorial-system-madesimpler-incase(iii)}
When condition (iii) of theorem \ref{t.mutation-stop-modifiedv4} above holds, the description of the system $\mathfrak{B}_{A,f}\mathfrak{Y}=\mathfrak{e}$ (or of the system $\mathfrak{B}_{A,f}^{tr}\mathfrak{Y}^{tr}=\mathfrak{e}^{tr}$ as the case may be) of linear equations as mentioned in subsection \ref{ss.description-of-anewsys} above can be made simpler. We will now give this simpler description by breaking up the matrix $\mathfrak{B}_{A,f}$ (or the matrix $\mathfrak{B}_{A,f}^{tr}$ as the case may be) into two parts: One corresponding to the columns of type A and the other corresponding to the columns of type B. We will henceforth call the part of the matrix corresponding to the type A elements as \textit{Part A of the matrix} $\mathfrak{B}_{A,f}$ (or of the matrix $\mathfrak{B}_{A,f}^{tr}$ as the case may be) and the rest as \textit{Part B}. Since we have assumed that condition (ii) of theorem \ref{t.mutation-stop-modifiedv4} does not hold, it follows that $(-2)_{max,A}<1_{min,A}$ in all the cases to be considered below.\\
\begin{center}
\underline{Case I: When $(-2)_{max,A}\geq 1$}
\end{center}
\noindent \underline{Subcase(I.1)}: When $(-3)_{max,A}\geq 1_{min,A}+2_{min,A}-1$.\\
\noindent Part A of the matrix $\mathfrak{B}_{A,f}$ is given by the matrix in table $1$.\\
\begin{table}[hbt]
\begin{center}
\caption{Part A of the matrix $\mathfrak{B}_{A,f}$} 
\begin{tabular}{|c|c|c|c|c|c|c|c|c|c|c|c|}
\hline
$\cdot$ & $\cdots$ & $\cdot$ & $\cdot$ & $\cdot$ & $\cdot$ & $A[\frac{-3}{1}]^{1_{min,A}}$ & $\cdot$ & $\cdot$ & $\cdots$ & $\cdot$ & $\leftarrow columns/rows\downarrow$ \\
\hline
 &  &  &  &  &  &  &  &  &  & $1$ & $\cdot$ \\
\hline
 &  &  &  &  &  &  &  &  & $\diagup$ &  & $\vdots$ \\
\hline
 &  &  &  &  &  &  &  & $1$ &  &  & $\cdot$ \\
\hline
 &  &  &  &  &  &  & $1$ &  &  &  & $\cdot$ \\
\hline
 &  &  &  &  &  & $-1$ &  &  &  &  & $A[\frac{-3}{1}]^{(1_{min,A}-1)}[\frac{-3}{2}]$ \\
\hline
 &  &  &  &  & $1$ & $1$ &  &  &  &  & $\cdot$ \\
\hline
 &  &  &  & $-1$ & $-2$ & $-1$ &  &  &  &  & $\cdot$ \\
\hline
 &  &  & $1$ & $3$ & $3$ & $1$ &  &  &  &  & $\cdot$ \\
\hline
 &  & $-1$ & $-4$ & $-6$ & $-4$ & $-1$ &  &  &  &  & $\cdot$ \\
\hline
 & $\vdots$ & $\cdot$ & $\cdot$ & $\cdot$ & $\cdot$ & $Alternate$ &  &  &  &  & $\vdots$ \\
\hline
 & $\cdot$ & $\cdot$ & $\cdot$ & $\cdot$ & $\cdot$ & $string\ of$ &  &  &  &  & $\vdots$ \\
\hline
$\cdot$ & $\cdot$ & $\cdot$ & $\cdot$ & $\cdot$ & $\cdot$ & $-1\ and\ 1$ &  &  &  &  & $\cdot$ \\
\hline
$\cdot$ & $\cdot$ & $\cdot$ & $\cdot$ & $\cdot$ & $\cdot$ & $without$ &  &  &  &  & $A[\frac{-3}{2}]$ \\
\hline
$\cdot$ & $\cdot$ & $\cdot$ & $\cdot$ & $\cdot$ & $\cdot$ & $gap$ &  &  &  &  & $A[\frac{-2}{1}]$ \\
\hline
\end{tabular}
\end{center}
\end{table}
\textit{Illustration of the matrix in table 1}: The columns of this matrix are indexed by the elements $A[\frac{-2}{1}]^{(-2)_{max,A}}[\frac{-3}{1}]^{1_{min,A}-(-2)_{max,A}}$ , $\cdots$ , $A[\frac{(-2)}{(1)}]^{2}[\frac{(-3)}{(1)}]^{1_{min,A}-2}$ , $A[\frac{(-2)}{(1)}]^{1}[\frac{(-3)}{(1)}]^{1_{min,A}-1}$ , $A[\frac{(-3)}{(1)}]^{1_{min,A}}$, $A[\frac{(-3)}{(1)}]^{1_{min,A}}[\frac{(-3)}{(2)}]$ , $A[\frac{(-3)}{(1)}]^{1_{min,A}}[\frac{(-3)}{(2)}]^{2}$ , $\cdots$ , $A[\frac{(-3)}{(1)}]^{1_{min,A}}[\frac{(-3)}{(2)}]^{2_{min,A}-1}$ from left to right in the `increasing' order induced by the order $\precsim$ on $\mathcal{L}(A,f)$. The rows of this matrix are indexed by the elements $A[\frac{-3}{1}]^{(1_{min,A})}[\frac{-3}{2}]^{2_{min,A}-1}$, $\cdots$, $A[\frac{-3}{1}]^{1_{min,A}}[\frac{-3}{2}]^{2}$, $A[\frac{-3}{1}]^{1_{min,A}}[\frac{-3}{2}]$, $A[\frac{-3}{1}]^{1_{min,A}-1}[\frac{-3}{2}]$, $A[\frac{-3}{1}]^{1_{min,A}-2}[\frac{-3}{2}]$, $\cdots$, $A[\frac{-3}{1}][\frac{-3}{2}]$, $A[\frac{-3}{2}]$, $A[\frac{-2}{1}]$ of the set $\mathcal{R}(A,f)$ which are arranged from top to bottom in the `decreasing' order induced by the order $\precsim$ on $\mathcal{L}(A,f)$. The blank spaces indicate $0$ entries. The row below the one containing the negative binomial coefficients of $4$ contains the positive binomial coefficients of $5$, the row below it contains the negative binomial coefficients of $6$, and so on. This string of signed binomial coefficients have one end at the column corresponding to the element $A[\frac{-3}{1}]^{1_{min,A}}$ and the other end is towards the left boundary of this table. The number of columns to the left of the column marked $A[\frac{-3}{1}]^{1_{min,A}}$ may or may not be sufficient for this string in a given row. But if there are sufficient such columns, then consider the cell where this string ends. For all existing cells towards the left of this cell, the corresponding entry is $0$.\\

\noindent Part B of the matrix $\mathfrak{B}_{A,f}$ is given by the matrix in table $2$.\\
\begin{table}[hbt]
\begin{center}
\caption{Part B of the matrix $\mathfrak{B}_{A,f}$} 
\begin{tabular}{|c|c|c|c|c|c|}
\hline
$\star_0$ & $\star_1$ & $\star_2$ & $\cdots$ & $\star_{(1_{min,A}-1)}$ & $\leftarrow columns/rows\downarrow$ \\
\hline
$\Diamond$ & $\Diamond$ & $\Diamond$ & $\Diamond$ & $\Diamond$ & $\blacklozenge_{(2_{min,A}-1)}$ \\
\hline
$\vdots$ & $\vdots$ & $\vdots$ & $\vdots$ & $\vdots$ & $\vdots$ \\
\hline
$\vdots$ & $\vdots$ & $\vdots$ & $\vdots$ & $\vdots$ & $\vdots$ \\
\hline
$\Diamond$ & $\Diamond$ & $\Diamond$ & $\Diamond$ & $\Diamond$ & $\blacklozenge_1$ \\
\hline
\hline
$\alpha_{(1_{min,A}-1)}$ & $\alpha_{(1_{min,A}-2)}$ & $\alpha_{(1_{min,A}-3)}$ & $\cdot$ & $\alpha_0$ & $\bigstar_{(1_{min,A}-1)}$ \\
\hline
$\cdot$ & $\cdot$ & $\cdot$ & $\cdot$ & $0$ & $\cdot$ \\
$\vdots$ & $\vdots$ & $\vdots$ & $\vdots$ & $\vdots$ & $\vdots$ \\
\hline
$\alpha_2$ & $\alpha_1$ & $\alpha_0$ & $\cdot$ & $0$ & $\cdot$ \\
\hline
$\alpha_1$ & $\alpha_0$ & $0$ & $\cdots$ & $0$ & $\bigstar_1$ \\
\hline
$\alpha_0$ & $0$ & $0$ & $\cdots$ & $0$ & $\bigstar_0$ \\
\hline
$0$ & $0$ & $0$ & $\cdots$ & $0$ & $A[\frac{-2}{1}]$ \\
\hline
\end{tabular}
\end{center}
\end{table}
\textit{Illustration of the matrix in table 2}: In this table,
\begin{center}
$\alpha_i:=(-1)^{(2_{min,A}+1)}\ ^{(2_{min,A}-1)}C_i$ for each $0\leq i\leq 2_{min,A}-1$,\\
$\star_i:=A[\frac{-3}{2}]^{2_{min,A}}[\frac{-3}{1}]^i$ for all $0\leq i\leq 1_{min,A}-1$, \\
$\bigstar_i:=A[\frac{-3}{1}]^{i}[\frac{-3}{2}]$ for all $0\leq i\leq 1_{min,A}-1$, \\
$\blacklozenge_i:=A[\frac{-3}{1}]^{1_{min,A}}[\frac{-3}{2}]^i$ for all $1\leq i\leq 2_{min,A}-1$ and\\
$\Diamond$ denotes some real entry, which does not\\
represent the same entry everywhere.
\end{center}
The columns of this matrix are indexed by the elements $A[\frac{-3}{2}]^{2_{min,A}}[\frac{-3}{1}]^i$ where $0\leq i\leq 1_{min,A}-1$, from left to right in the `increasing' order induced by the order $\precsim$ on $\mathcal{L}(A,f)$. The rows of this matrix are indexed similarly as in table $1$. In the lowest row, all the entries are $0$. Each column contains a string (without gap) of elements $\alpha_0,\alpha_1,\ldots,$ where $\alpha_0$ is the bottom-most, above it $\alpha_1$, $\ldots$, and so on. This string truncates at the double-lined partition indicated in the table. In the left-most column, this string begins at the $2$-nd row from the bottom, where by the phrase ` this string begins' we mean that $\alpha_0$ is the entry of the indicated position. In the next ($2$-nd from left to right) column, this string begins at the $3$-rd row from the bottom, and so on, till the last column where this string begins at the row corresponding to the index $\bigstar_{(1_{min,A}-1)}$. In each column, the entries below the element $\alpha_0$ are all $0$. It is not necessary to know the entries of this table which are denoted by the symbol $\Diamond$, the reason behind this is explained in remark \ref{r.no-need-toknow-entries}. \\ 

\noindent \underline{Subcase(I.2)}: When $max\{1_{min,A},2_{min,A}\}\leq (-3)_{max,A}<1_{min,A}+2_{min,A}-1$.\\
For getting Part A of the matrix $\mathfrak{B}_{A,f}$, remove all those rows and columns from the matrix in table $1$ which are indexed by elements containing negative powers of any of the underlying variables. The resulting matrix is Part A of the matrix $\mathfrak{B}_{A,f}$. Apply a similar process on the matrix in table $2$ to get the Part B of the matrix $\mathfrak{B}_{A,f}$.\\

\noindent \underline{Subcase(I.3)}: When $(-3)_{max,A}<max\{1_{min,A},2_{min,A}\}$.\\

\noindent Possibility (a): When $(-3)_{max,A}<2_{min,A}$ and $(-3)_{max,A}\geq 1_{min,A}$.\\
The matrix $\mathfrak{B}_{A,f}$ consists only of Part A, there is no Part B. The structure of Part A of the matrix $\mathfrak{B}_{A,f}$ is similar that in subcase (I.2).\\

\noindent Possibility (b): When $(-3)_{max,A}<2_{min,A}$ and $(-3)_{max,A}<1_{min,A}\leq (-2)_{max,A}+(-3)_{max,A}$.\\
The matrix $\mathfrak{B}_{A,f}$ consists only of Part A, there is no Part B. To get Part A of the matrix $\mathfrak{B}_{A,f}$, remove all those rows from the matrix in table $1$ which are indexed by elements containing negative powers of any of the underlying variables. And keep only those columns of the matrix in table $1$ which are indexed by the elements $A[\frac{-2}{1}]^{(-2)_{max,A}}[\frac{-3}{1}]^{1_{min,A}-(-2)_{max,A}}$ , $\cdots$ ,$A[\frac{-2}{1}]^{1_{min,A}-(-3)_{max,A}}[\frac{-3}{1}]^{(-3)_{max,A}}$ from left to right in the `increasing' order induced by the order $\precsim$ on $\mathcal{L}(A,f)$. The resulting matrix is Part A of the matrix $\mathfrak{B}_{A,f}$. \\

\noindent Possibility (c): When $(-3)_{max,A}\geq 2_{min,A}$. \\
In this situation, we must have $(-3)_{max,A}<1_{min,A}\leq (-2)_{max,A}+(-3)_{max,A}$. Both Part A and Part B of the matrix $\mathfrak{B}_{A,f}$ will exist non trivially. The structure of Part A of the matrix $\mathfrak{B}_{A,f}$ is similar to that in possibility (b). The structure of Part B of the matrix $\mathfrak{B}_{A,f}$ is similar to that in subcase (I.2).
\begin{center}
\underline{Case II: When $(-2)_{max,A}=0$}
\end{center}
I provide here a detailed description of the equivalent system $\mathfrak{B}_{A,f}^{tr}\mathfrak{Y}^{tr}=\mathfrak{e}^{tr}$ of linear equations.\\
\noindent \underline{Subcase(II.1)}: When $(-3)_{max,A}\geq 1_{min,A}+2_{min,A}-1$.\\
\noindent Part A of the matrix $\mathfrak{B}_{A,f}^{tr}$ is given by the matrix in table $3$.\\
\begin{table}[hbt]
\begin{center}
\caption{Part A of the matrix $\mathfrak{B}_{A,f}^{tr}$} 
\begin{tabular}{|c|c|c|c|c|c|}
\hline
$A[\frac{-3}{1}]^{1_{min,A}}$ & $\cdot$ & $\cdot$ & $\cdots$ & $\cdot$ & $\leftarrow columns/rows\downarrow$ \\
\hline
  &  &  &  & $1$ & $\cdot$ \\
\hline
  &  &  & $\diagup$ &  & $\vdots$ \\
\hline
  &  & $1$ &  &  & $\cdot$ \\
\hline
  & $1$ &  &  &  & $\cdot$ \\
\hline
$-1$ &  &  &  &  & $A[\frac{-3}{1}]^{(1_{min,A}-1)}[\frac{-3}{2}]$ \\
\hline
$1$ &  &  &  &  & $\cdot$ \\
\hline
$-1$ &  &  &  &  & $\cdot$ \\
\hline
$1$ &  &  &  &  & $\cdot$ \\
\hline
$-1$ &  &  &  &  & $\cdot$ \\
\hline
$Alternate$ &  &  &  &  & $\vdots$ \\
\hline
$string\ of$ &  &  &  &  & $\vdots$ \\
\hline
$-1\ and\ 1$ &  &  &  &  & $\cdot$ \\
\hline
$without$ &  &  &  &  & $A[\frac{-3}{2}]$ \\
\hline
$gap$ &  &  &  &  & $A[\frac{-3}{1}]$ \\
\hline
\end{tabular}
\end{center}
\end{table}
\textit{Illustration of the matrix in table 3}: The columns of this matrix are indexed by the elements $A[\frac{(-3)}{(1)}]^{1_{min,A}}$, $A[\frac{(-3)}{(1)}]^{1_{min,A}}[\frac{(-3)}{(2)}]$ , $A[\frac{(-3)}{(1)}]^{1_{min,A}}[\frac{(-3)}{(2)}]^{2}$ , $\cdots$ , $A[\frac{(-3)}{(1)}]^{1_{min,A}}[\frac{(-3)}{(2)}]^{2_{min,A}-1}$ from left to right in the `increasing' order induced by the order $\precsim$ on $\mathcal{L}(A,f)$. The rows of this matrix are indexed by the elements $A[\frac{-3}{1}]^{(1_{min,A})}[\frac{-3}{2}]^{2_{min,A}-1}$, $\cdots$, $A[\frac{-3}{1}]^{1_{min,A}}[\frac{-3}{2}]^{2}$, $A[\frac{-3}{1}]^{1_{min,A}}[\frac{-3}{2}]$, $A[\frac{-3}{1}]^{1_{min,A}-1}[\frac{-3}{2}]$, $A[\frac{-3}{1}]^{1_{min,A}-2}[\frac{-3}{2}]$, $\cdots$, $A[\frac{-3}{1}][\frac{-3}{2}]$, $A[\frac{-3}{2}]$, $A[\frac{-3}{1}]$ of the set $\mathcal{R}(A,f)$ which are arranged from top to bottom in the `decreasing' order induced by the order $\precsim$ on $\mathcal{L}(A,f)$. The blank spaces indicate $0$ entries. \\

\noindent Part B of the matrix $\mathfrak{B}_{A,f}^{tr}$ is of the same form as Part B of the matrix $\mathfrak{B}_{A,f}$ when $(-2)_{max,A}\geq 1$ and $(-3)_{max,A}\geq 1_{min,A}+2_{min,A}-1$, the only difference being that the last row here is indexed by $A[\frac{-3}{1}]$ instead of $A[\frac{-2}{1}]$.\\

\noindent \underline{Subcase(II.2)}: When $max\{1_{min,A},2_{min,A}\}\leq (-3)_{max,A}<1_{min,A}+2_{min,A}-1$.\\
For getting Part A of the matrix $\mathfrak{B}_{A,f}^{tr}$, remove all those rows and columns from the matrix in table $3$ which are indexed by elements containing negative powers of any of the underlying variables. The resulting matrix is Part A of the matrix $\mathfrak{B}_{A,f}^{tr}$. For getting Part B of the matrix $\mathfrak{B}_{A,f}^{tr}$, consider Part B of the matrix $\mathfrak{B}_{A,f}^{tr}$ of subcase (II.1) and apply a similar process on it (as we did for Part A).\\

\noindent \underline{Subcase(II.3)}: When $(-3)_{max,A}<max\{1_{min,A},2_{min,A}\}$.\\

\noindent Possibility (a): When $(-3)_{max,A}<2_{min,A}$ but $(-3)_{max,A}\geq 1_{min,A}$. \\
The matrix $\mathfrak{B}_{A,f}^{tr}$ consists only of Part A, there is no Part B. The structure of Part A of the matrix $\mathfrak{B}_{A,f}^{tr}$ is similar to that in subcase(II.2). \\

\noindent Possibility (b): When $(-3)_{max,A}<2_{min,A}$ and $1_{min,A}$ both. \\
Neither Part A nor Part B of the matrix $\mathfrak{B}_{A,f}^{tr}$ exists.\\

\noindent Possibility (c): When $(-3)_{max,A}\geq 2_{min,A}$ and $(-3)_{max,A}<1_{min,A}$. \\
The matrix $\mathfrak{B}_{A,f}^{tr}$ consists only of Part B, there is no Part A. The structure of Part B of the matrix $\mathfrak{B}_{A,f}^{tr}$ is similar to that in subcase(II.2). \\
\section{Solvability of the final system of equations}\label{s.solvability-of-thenewsys}
Recall from theorem \ref{t.mutation-stop-modifiedv4} above that the main problem has now reduced to checking the solvability of certain systems of linear equations. We will now transform the system $\mathfrak{B}_{A,f}\mathfrak{Y}=\mathfrak{e}$ (or $\mathfrak{B}_{A,f}^{tr}\mathfrak{Y}^{tr}=\mathfrak{e}^{tr}$ as the case may be) of equations into an equivalent and simpler system $\mathfrak{B}_{A,f}^{eq}\mathfrak{Y}^{eq}=\mathfrak{e}^{eq}$ of linear equations which will help us decide about the solvability of the system $\mathfrak{B}_{A,f}\mathfrak{Y}=\mathfrak{e}$ or $\mathfrak{B}_{A,f}^{tr}\mathfrak{Y}^{tr}=\mathfrak{e}^{tr}$. 

In all the tables that we come across hereafter, the lines along the south-west$\leftrightarrow$north-east direction indicate the continuation of the same entry along that direction. That is, the entries in any two adjacent cells along that line are the same.
\begin{center}
\underline{Case I: When $(-2)_{max,A}\geq 1$}
\end{center}
\noindent \underline{Subcase(I.1)}: When $(-3)_{max,A}\geq 1_{min,A}+2_{min,A}-1$.\\
Let us apply the following row operations on the matrix $\mathfrak{B}_{A,f}$:\\
\noindent Stage 1: For any row that lies strictly below the row indexed by $A[\frac{-3}{1}]^{(1_{min,A}-1)}[\frac{-3}{2}]$, replace it with `the row just above it $+$ itself'.\\
\noindent Stage 2: For any row that lies strictly below the row indexed by $A[\frac{-3}{1}]^{(1_{min,A}-2)}[\frac{-3}{2}]$, replace it with `the row just above it $+$ itself'.\\
$\ldots\ldots$ and so on. Proceed similarly upto Stage $(-2)_{max,A}+1$. This will result into the equivalent system $\mathfrak{B}_{A,f}^{eq}\mathfrak{Y}^{eq}=\mathfrak{e}^{eq}$ of linear equations where $\mathfrak{Y}^{eq}=\mathfrak{Y}$ and $\mathfrak{e}^{eq}=\mathfrak{e}$. Part A of the matrix $\mathfrak{B}_{A,f}^{eq}$ has the property that the submatrix of this matrix formed by the last $1_{min,A}-(-2)_{max,A}$ many rows of it is the zero matrix. The submatrix of Part A of the matrix $\mathfrak{B}_{A,f}^{eq}$ formed by removing the last $1_{min,A}-(-2)_{max,A}$ many rows of it is a matrix having $1$ s or $-1$ s on the antidiagonal and all other entries $0$. Therefore it is not needed to know Part B of the matrix $\mathfrak{B}_{A,f}^{eq}$ fully, see remark \ref{r.no-need-toknow-entries} for a reasoning. It is enough to know the submatrix of Part B of the matrix $\mathfrak{B}_{A,f}^{eq}$ which is formed by the last $1_{min,A}-(-2)_{max,A}$ many rows of it. Let us denote this submatrix by $\mathcal{C}_{A,f}$, table $4$ provides a description of this matrix.\\   
\begin{table}[hbt]
\begin{center}
\caption{The matrix $\mathfrak{C}_{A,f}$ of subcase(I.1)} 
\begin{tabular}{|c|c|c|c|c|c|c|c|}
\hline
$\star_0$ & $\star_1$ & $\cdots$ & $\star_{(-2)_{max,A}+2}$ & $\cdots$ & $\star_{(1_{min,A}-2)}$ & $\star_{(1_{min,A}-1)}$ & $\leftarrow cols/rows\downarrow$ \\
\hline
$\beta_{(1_{min,A}-1)}$ & $\beta_{(1_{min,A}-2)}$ & $\cdots$ & $\cdots$ & $\cdots$ & $\beta_1$ & $\beta_0$ & $\Join$ \\
\hline
$\beta_{(1_{min,A}-2)}$ & $\diagup$ & $\diagup$ & $\diagup$ & $\beta_1$ & $\beta_0$ & $0$ & $\cdot$ \\
\hline
$\diagup$ & $\diagup$ & $\diagup$ & $\diagup$ & $\diagup$ & $0$ & $0$ & $\vdots$ \\
\hline
$\cdots$ & $\beta_2$ & $\beta_1$ & $\beta_0$ & $\cdots$ & $0$ & $0$ & $A[\frac{-2}{1}]$ \\
\hline
\end{tabular}
\end{center}
\end{table}
\textit{Illustration of the matrix in table 4}: In this table,
\begin{center}
$\bigstar_i:=A[\frac{-3}{1}]^{i}[\frac{-3}{2}]$ for all $0\leq i\leq 1_{min,A}-1$, \\
$\Join:=\bigstar_{1_{min,A}-(-2)_{max,A}-2}$,\\
$\alpha_i:=(-1)^{(2_{min,A}+1)}\ ^{(2_{min,A}-1)}C_i$ for each $0\leq i\leq 2_{min,A}-1$,\\
$\ ^{((-2)_{max,A}+1)}C_r=0$ if $r>(-2)_{max,A}+1$, \\
$\beta_j:=(\alpha_0,\ldots,\alpha_{j-1},\alpha_{j})\rtimes(\delta_{j},\ldots,\delta_1,\delta_0)$ for each $0\leq j\leq 1_{min,A}-1$\\
where for each $0\leq j\leq 1_{min,A}-1$, $\delta_{l}:=\ ^{((-2)_{max,A}+1)}C_l$ if $0\leq l\leq j$\\
and $\rtimes$ denotes the `standard inner product' of the two vectors\\
$(\alpha_0,\ldots,\alpha_{j-1},\alpha_{j})$ and $(\delta_{j},\ldots,\delta_1,\delta_0)$.\\
$\star_i:=A[\frac{-3}{2}]^{2_{min,A}}[\frac{-3}{1}]^i$ for all $0\leq i\leq 1_{min,A}-1$ and\\
In any column containing $\beta_0$, all elements below $\beta_0$ are zero.
\end{center}
The system $\mathfrak{B}_{A,f}^{eq}\mathfrak{Y}^{eq}=\mathfrak{e}^{eq}$ of linear equations is solvable if and only if the rank of the matrix $\mathcal{C}_{A,f}$ equals the rank of the matrix $[\mathcal{C}_{A,f}\ \mathcal{E}]$ where $\mathcal{E}$ is the vector $(0,\ldots,0,1)^t$. But these two ranks are equal since $\beta_0\neq 0$, hence the system $\mathfrak{B}_{A,f}^{eq}\mathfrak{Y}^{eq}=\mathfrak{e}^{eq}$ is solvable.\\

\noindent \underline{Subcase(I.2)}:  When $max\{1_{min,A},2_{min,A}\}\leq (-3)_{max,A}<1_{min,A}+2_{min,A}-1$.\\
The equivalent system $\mathfrak{B}_{A,f}^{eq}\mathfrak{Y}^{eq}=\mathfrak{e}^{eq}$ is obtained from the original system $\mathfrak{B}_{A,f}\mathfrak{Y}=\mathfrak{e}$ in exactly the same way as in subcase (I.1). Due to a similar reason as in subcase (I.1), it is enough to look at the submatrix of Part B of the matrix $\mathfrak{B}_{A,f}^{eq}$ formed by the last $1_{min,A}-(-2)_{max,A}$ many rows of it. Let us denote this submatrix by $\mathcal{C}_{A,f}$, table $5$ provides a description of this matrix.\\
\begin{table}[hbt]
\begin{center}
\caption{The matrix $\mathfrak{C}_{A,f}$ of subcase(I.2)} 
\begin{tabular}{|c|c|c|c|c|c|c|c|}
\hline
$\star_0$ & $\star_1$ & $\cdots$ & $\cdots$ & $\cdots$ & $\star_{k_2-1}$ & $\star_{k_2}$ & $\leftarrow cols/rows\downarrow$ \\
\hline
$\beta_{(1_{min,A}-1)}$ & $\beta_{(1_{min,A}-2)}$ & $\cdots$ & $\cdots$ & $\cdots$ & $\beta_{(1_{min,A}-k_2)}$ & $\beta_{(1_{min,A}-k_2-1)}$ & $\Join$ \\
\hline
$\beta_{(1_{min,A}-2)}$ & $\diagup$ & $\diagup$ & $\diagup$ & $\diagup$ & $\diagup$ & $\beta_{(1_{min,A}-k_2-2)}$ & $\cdot$ \\
\hline
$\diagup$ & $\diagup$ & $\diagup$ & $\diagup$ & $\diagup$ & $\diagup$ & $\vdots$ & $\vdots$ \\
\hline
$\diagup$ & $\diagup$ & $\diagup$ & $\diagup$ & $\diagup$ & $\diagup$ & $\cdot$ & $A[\frac{-2}{1}]$ \\
\hline
\end{tabular}
\end{center}
\end{table}
\textit{Illustration of the matrix in table 5}: In this table, $(-3)_{max,A}=2_{min,A}+k_2$ and all other notation remains the same as in table $4$. In each column, there is a string of $\beta_j$-s which begins from the topmost row. This string has the property that if the entry in any fixed row at any column is $\beta_j$, then the entry in the row just below that row at the same column is $\beta_{j-1}$ if $j\geq 1$ and $0$ otherwise.\\
The system $\mathfrak{B}_{A,f}^{eq}\mathfrak{Y}^{eq}=\mathfrak{e}^{eq}$ of linear equations is solvable if and only if the rank of the matrix $\mathcal{C}_{A,f}$ equals the rank of the matrix $[\mathcal{C}_{A,f}\ \mathcal{E}]$ where $\mathcal{E}$ is the vector $(0,\ldots,0,1)^t$. But unlike in subcase (I.1), it is not clear here whether or not the ranks of these two matrices are equal. Hence one needs to compute the ranks of such matrices to determine whether or not the system $\mathfrak{B}_{A,f}^{eq}\mathfrak{Y}^{eq}=\mathfrak{e}^{eq}$ is solvable.\\

\noindent \underline{Subcase(I.3)}: When $(-3)_{max,A}<max\{1_{min,A},2_{min,A}\}$.\\

\noindent Possibility (a): When $(-3)_{max,A}<2_{min,A}$ and $(-3)_{max,A}\geq 1_{min,A}$.\\
Since the matrix $\mathfrak{B}_{A,f}$ consists only of Part A (there is no Part B), the same is true for the matrix $\mathfrak{B}_{A,f}^{eq}$. The structure of Part A of the matrix $\mathfrak{B}_{A,f}^{eq}$ is similar that in subcase (I.2). It is clear from the structure of this matrix that the system $\mathfrak{B}_{A,f}^{eq}\mathfrak{Y}^{eq}=\mathfrak{e}^{eq}$ of linear equations is solvable if and only if $1_{min,A}\leq (-2)_{max,A}$.\\

\noindent Possibility (b): When $(-3)_{max,A}<2_{min,A}$ and $(-3)_{max,A}<1_{min,A}\leq (-2)_{max,A}+(-3)_{max,A}$.\\
Since the matrix $\mathfrak{B}_{A,f}$ consists only of Part A (there is no Part B), the same is true for the matrix $\mathfrak{B}_{A,f}^{eq}$. Suppose $(-3)_{max,A}=1_{min,A}-k$. Clearly then $k\leq (-2)_{max,A}$. Let us apply the following row operations on the matrix $\mathfrak{B}_{A,f}$:\\
\noindent Stage 1: For any row that lies strictly below the row indexed by $A[\frac{-3}{1}]^{(1_{min,A}-k-1)}[\frac{-3}{2}]$, replace it with `the row just above it $+$ itself'.\\
\noindent Stage 2: Perform the same activities as in Stage 1 $(k+1)$-many times.\\
\noindent Stage 3: For any row that lies strictly below the row indexed by $A[\frac{-3}{1}]^{(1_{min,A}-k-2)}[\frac{-3}{2}]$, replace it with `the row just above it $+$ itself'.\\
\noindent Stage 4: For any row that lies strictly below the row indexed by $A[\frac{-3}{1}]^{(1_{min,A}-k-3)}[\frac{-3}{2}]$, replace it with `the row just above it $+$ itself'.\\
$\ldots\ldots$ and so on. Proceed similarly (as in Stage 3) upto Stage $(-2)_{max,A}-k+2$. It should be noted here that the row operations performed upto stage 2 are different than those performed in stage 3 and afterwards. This will result into the equivalent system $\mathfrak{B}_{A,f}^{eq}\mathfrak{Y}^{eq}=\mathfrak{e}^{eq}$ of linear equations where $\mathfrak{Y}^{eq}=\mathfrak{Y}$ and $\mathfrak{e}^{eq}=\mathfrak{e}$. Part A of the matrix $\mathfrak{B}_{A,f}^{eq}$ has the property that the submatrix of it formed by the last $1_{min,A}-(-2)_{max,A}$ many rows of it is the zero matrix. The submatrix of Part A of the matrix $\mathfrak{B}_{A,f}^{eq}$ formed by removing the last $1_{min,A}-(-2)_{max,A}$ many rows of it is a matrix having $1$ s or $-1$ s on the antidiagonal and all other entries $0$. It is now clear that the system $\mathfrak{B}_{A,f}^{eq}\mathfrak{Y}^{eq}=\mathfrak{e}^{eq}$ of linear equations is solvable if and only if $1_{min,A}\leq (-2)_{max,A}$.\\

\noindent Possibility (c): When $(-3)_{max,A}\geq 2_{min,A}$. \\
In this situation, we must have $(-3)_{max,A}<1_{min,A}\leq (-2)_{max,A}+(-3)_{max,A}$. Unlike in the possibilities (a) and (b) above, both Part A and Part B of the matrix $\mathfrak{B}_{A,f}^{eq}$ exist non trivially. The structure of Part A of the matrix $\mathfrak{B}_{A,f}^{eq}$ is similar to that in possibility (b) of subcase (I.3). Let us now describe the structure of Part B of the matrix $\mathfrak{B}_{A,f}^{eq}$. Let $k$ and $k'$ be such that $(-3)_{max,A}=1_{min,A}-k$ and $(-3)_{max,A}=2_{min,A}+k'$. Due to a similar reason as in subcase (I.1), it is enough to look at the submatrix of Part B of the matrix $\mathfrak{B}_{A,f}^{eq}$ formed by the last $1_{min,A}-(-2)_{max,A}$ many rows of it. Let us denote this submatrix by $\mathcal{C}_{A,f}$, table $6$ provides a description of this matrix.\\
\begin{table}[hbt]
\begin{center}
\caption{The matrix $\mathfrak{C}_{A,f}$ of subcase(I.3,Possibility (c))} 
\begin{tabular}{|c|c|c|c|c|c|c|c|}
\hline
$\star_0$ & $\star_1$ & $\cdots$ & $\cdots$ & $\cdots$ & $\star_{k'-1}$ & $\star_{k'}$ & $\leftarrow cols/rows\downarrow$ \\
\hline
$\heartsuit_k^0$ & $\heartsuit_k^1$ & $\cdots$ & $\cdots$ & $\cdots$ & $\heartsuit_k^{k'-1}$ & $\heartsuit_k^{k'}$ & $\square_{(1_{min,A}-(-2)_{max,A}-2)}$ \\
\hline
$\vdots$ & $\vdots$ & $\vdots$ & $\vdots$ & $\vdots$ & $\vdots$ & $\vdots$ & $\vdots$ \\
\hline
$\heartsuit_2^0$ & $\heartsuit_2^1$ & $\cdots$ & $\cdots$ & $\cdots$ & $\heartsuit_2^{k'-1}$ & $\heartsuit_2^{k'}$ & $\square_{(1_{min,A}-(-2)_{max,A}-k)}$ \\
\hline
$\heartsuit_1^0$ & $\heartsuit_1^1$ & $\cdots$ & $\cdots$ & $\cdots$ & $\heartsuit_1^{k'-1}$ & $\heartsuit_1^{k'}$ & $\square_{(1_{min,A}-(-2)_{max,A}-k-1)}$ \\
\hline
$\blacksquare_0$ & $\blacksquare_1$ & $\cdots$ & $\cdots$ & $\cdots$ & $\blacksquare_{k'-1}$ & $\blacksquare_{k'}$ & $\square_{(1_{min,A}-(-2)_{max,A}-k-2)}$ \\
\hline
$\blacksquare_1$ & $\blacksquare_2$ & $\cdots$ & $\cdots$ & $\cdots$ & $\blacksquare_{k'}$ & $\blacksquare_{k'+1}$ & $\square_{(1_{min,A}-(-2)_{max,A}-k-3)}$ \\
\hline
$\diagup$ & $\diagup$ & $\diagup$ & $\diagup$ & $\diagup$ & $\diagup$ & $\vdots$ & $\vdots$ \\
\hline
$\diagup$ & $\diagup$ & $\diagup$ & $\diagup$ & $\diagup$ & $\diagup$ & $\vdots$ & $\square_0$ \\
\hline
$\diagup$ & $\diagup$ & $\diagup$ & $\diagup$ & $\diagup$ & $\diagup$ & $\cdot$ & $A[\frac{-2}{1}]$ \\
\hline
\end{tabular}
\end{center}
\end{table}
\linebreak
\textit{Illustration of the matrix in table 6}: In this table,
\begin{center}
$\square_i:=A[\frac{-3}{1}]^{i}[\frac{-3}{2}]$ for all $0\leq i\leq 1_{min,A}-k-1$,\\
$\star_i:=A[\frac{-3}{2}]^{2_{min,A}}[\frac{-3}{1}]^i$ for all $0\leq i\leq 1_{min,A}-1$,\\
$\alpha_i:=(-1)^{(2_{min,A}+1)}\ ^{(2_{min,A}-1)}C_i$ for each $0\leq i\leq 2_{min,A}-1$,\\
$\delta_{l}:=\ ^{((-2)_{max,A}+1)}C_l$ for any $0\leq l\leq (-2)_{max,A}+1$\\
and $\delta_l:=0$ if $l>(-2)_{max,A}+1$.\\
$\blacksquare_j:=(\alpha_0,\ldots,\alpha_{1_{min,A}-(k+1)-j-1},\alpha_{1_{min,A}-(k+1)-j})\rtimes(\delta_{1_{min,A}-(k+1)-j},\ldots,\delta_1,\delta_0)$\\
for each $0\leq j\leq 1_{min,A}-(k+1)$ and $\blacksquare_j:=0$ if $j>1_{min,A}-(k+1)$\\
where $\rtimes$ denotes the `standard inner product' of the two vectors\\
$(\alpha_0,\ldots,\alpha_{1_{min,A}-(k+1)-j-1},\alpha_{1_{min,A}-(k+1)-j})$ and $(\delta_{1_{min,A}-(k+1)-j},\ldots,\delta_1,\delta_0)$.\\
For any $1\leq i\leq k$ and $0\leq r\leq k'$,$\heartsuit_i^r:=$\\
$(\alpha_0,\ldots,\alpha_{1_{min,A}-(k+1)-r-1},\alpha_{1_{min,A}-(k+1)-r})\rtimes(\delta_{1_{min,A}-(k+1)+i-r},\ldots,\delta_{i+1},\delta_i)$.\\
\end{center}
In each column, the subscript of $j$ the elements $\blacksquare_j$ increases as we move from top to the bottom and the subscript $i$ of the elements $\heartsuit_i^r$ decreases as we move from top to the bottom. The superscript $r$ of the elements $\heartsuit_i^r$ ranges from $0$ to $k'$ and it increases as we move from left to right.\\

The system $\mathfrak{B}_{A,f}^{eq}\mathfrak{Y}^{eq}=\mathfrak{e}^{eq}$ of linear equations is solvable if and only if the rank of the matrix $\mathcal{C}_{A,f}$ equals the rank of the matrix $[\mathcal{C}_{A,f}\ \mathcal{E}]$ where $\mathcal{E}$ is the vector $(0,\ldots,0,1)^t$. Unlike in subcase (I.1), it is not clear here whether or not the ranks of these two matrices are equal. Hence one needs to compute the ranks of such matrices to determine whether or not the system $\mathfrak{B}_{A,f}^{eq}\mathfrak{Y}^{eq}=\mathfrak{e}^{eq}$ is solvable.\\
\begin{center}
\underline{Case II: When $(-2)_{max,A}=0$}
\end{center}
\noindent \underline{Subcase(II.1)}: When $(-3)_{max,A}\geq 1_{min,A}+2_{min,A}-1$.\\
Consider the system $\mathfrak{B}_{A,f}^{tr}\mathfrak{Y}^{tr}=\mathfrak{e}^{tr}$ of linear equations as mentioned in case II of subsection \ref{ss.new-combinatorial-system-madesimpler-incase(iii)} above. Let us apply the following row operations on the matrix $\mathfrak{B}_{A,f}^{tr}$:
\begin{center}
For any row that lies strictly below the row indexed by $A[\frac{-3}{1}]^{(1_{min,A}-1)}[\frac{-3}{2}]$, replace it with `the row just above it $+$ itself'.
\end{center}
This will result into the equivalent system $\mathfrak{B}_{A,f}^{eq}\mathfrak{Y}^{eq}=\mathfrak{e}^{eq}$ of linear equations where $\mathfrak{Y}^{eq}=\mathfrak{Y}^{tr}$ and $\mathfrak{e}^{eq}=\mathfrak{e}^{tr}$. Part A of the matrix $\mathfrak{B}_{A,f}^{eq}$ has the property that the submatrix of this matrix formed by the last $1_{min,A}$ many rows of it is the zero matrix. The submatrix of Part A of the matrix $\mathfrak{B}_{A,f}^{eq}$ formed by removing the last $1_{min,A}$ many rows of it is a matrix having $1$ s or $-1$ s on the antidiagonal and all other entries $0$. Therefore it is not needed to know Part B of the matrix $\mathfrak{B}_{A,f}^{eq}$ fully, see remark \ref{r.no-need-toknow-entries} for a reasoning. It is enough to know the submatrix of Part B of the matrix $\mathfrak{B}_{A,f}^{eq}$ which is formed by the last $1_{min,A}$ many rows of it. Let us denote this submatrix by $\mathcal{C}_{A,f}$, table $7$ provides a description of this matrix.\\   
\begin{table}[hbt]
\begin{center}
\caption{The matrix $\mathfrak{C}_{A,f}$ of subcase(II.1)} 
\begin{tabular}{|c|c|c|c|c|c|c|c|}
\hline
$\star_0$ & $\star_1$ & $\cdots$ & $\cdots$ & $\star_{(1_{min,A}-2)}$ & $\star_{(1_{min,A}-1)}$ & $\leftarrow cols/rows\downarrow$ \\
\hline
$\beta_{(1_{min,A}-1)}^{(0)}$ & $\beta_{(1_{min,A}-2)}^{(0)}$ & $\cdots$ & $\cdots$ & $\beta_1^{(0)}$ & $\beta_0^{(0)}$ & $\Join^{(0)}$ \\
\hline
$\beta_{(1_{min,A}-2)}^{(0)}$ & $\diagup$ & $\diagup$ & $\beta_1^{(0)}$ & $\beta_0^{(0)}$ & $0$ & $\cdot$ \\
\hline
$\diagup$ & $\diagup$ & $\diagup$ & $\diagup$ & $\diagup$ & $\vdots$ & $\vdots$ \\
\hline
$\diagup$ & $\diagup$ & $\diagup$ & $\diagup$ & $\cdot$ & $0$ & $A[\frac{-3}{1}]$ \\
\hline
\end{tabular}
\end{center}
\end{table}
\textit{Illustration of the matrix in table 7}: In this table,
\begin{center}
$\star_i:=A[\frac{-3}{2}]^{2_{min,A}}[\frac{-3}{1}]^i$ for all $0\leq i\leq 1_{min,A}-1$, \\
$\bigstar_i:=A[\frac{-3}{1}]^{i}[\frac{-3}{2}]$ for all $0\leq i\leq 1_{min,A}-1$, \\
$\Join^{(0)}:=\bigstar_{1_{min,A}-2}$,\\
$\alpha_i:=(-1)^{(2_{min,A}+1)}\ ^{(2_{min,A}-1)}C_i$ for each $0\leq i\leq 2_{min,A}-1$,\\
$\beta_j^{(0)}:=\alpha_{j-1}+\alpha_{j}$ for each $1\leq j\leq 1_{min,A}-1$ and $\beta_0^{(0)}:=\alpha_0$.\\
\end{center}
The system $\mathfrak{B}_{A,f}^{eq}\mathfrak{Y}^{eq}=\mathfrak{e}^{eq}$ of linear equations is solvable if and only if the rank of the matrix $\mathcal{C}_{A,f}$ equals the rank of the matrix $[\mathcal{C}_{A,f}\ \mathcal{E}]$ where $\mathcal{E}$ is the vector $(0,\ldots,0,1)^t$. But these two ranks are equal since $\beta_0^{(0)}\neq 0$, hence the system $\mathfrak{B}_{A,f}^{eq}\mathfrak{Y}^{eq}=\mathfrak{e}^{eq}$ is solvable.\\

\noindent \underline{Subcase(II.2)}: When $max\{1_{min,A},2_{min,A}\}\leq (-3)_{max,A}<1_{min,A}+2_{min,A}-1$.\\
The equivalent system $\mathfrak{B}_{A,f}^{eq}\mathfrak{Y}^{eq}=\mathfrak{e}^{eq}$ is obtained from the original system $\mathfrak{B}_{A,f}^{tr}\mathfrak{Y}^{tr}=\mathfrak{e}^{tr}$ in exactly the same way as in subcase (II.1). Due to a similar reason as in subcase (II.1), it is enough to look at the submatrix of Part B of the matrix $\mathfrak{B}_{A,f}^{eq}$ formed by the last $1_{min,A}$ many rows of it. Let us denote this submatrix by $\mathcal{C}_{A,f}$, table $8$ provides a description of this matrix.\\
\begin{table}[hbt]
\begin{center}
\caption{The matrix $\mathfrak{C}_{A,f}$ of subcase(II.2)} 
\begin{tabular}{|c|c|c|c|c|c|c|c|}
\hline
$\star_0$ & $\star_1$ & $\cdots$ & $\cdots$ & $\star_{k_2-1}$ & $\star_{k_2}$ & $\leftarrow cols/rows\downarrow$ \\
\hline
$\beta_{(1_{min,A}-1)}^{(0)}$ & $\beta_{(1_{min,A}-2)}^{(0)}$ & $\cdots$ & $\cdots$ & $\beta_{(1_{min,A}-k_2)}^{(0)}$ & $\beta_{(1_{min,A}-k_2-1)}^{(0)}$ & $\Join^{(0)}$ \\
\hline
$\beta_{(1_{min,A}-2)}^{(0)}$ & $\diagup$ & $\diagup$ & $\diagup$ & $\beta_{(1_{min,A}-k_2-1)}^{(0)}$ & $\beta_{(1_{min,A}-k_2-2)}^{(0)}$ & $\cdot$ \\
\hline
$\diagup$ & $\diagup$ & $\diagup$ & $\diagup$ & $\diagup$ & $\vdots$ & $\vdots$ \\
\hline
$\diagup$ & $\diagup$ & $\diagup$ & $\diagup$ & $\diagup$ & $\cdot$ & $A[\frac{-3}{1}]$ \\
\hline
\end{tabular}
\end{center}
\end{table}
\textit{Illustration of the matrix in table 8}: In this table,
\begin{center}
$k_2:=(-3)_{max,A}-2_{min,A}$,\\
$\star_i:=A[\frac{-3}{2}]^{2_{min,A}}[\frac{-3}{1}]^i$ for all $0\leq i\leq 1_{min,A}-1$, \\
$\bigstar_i:=A[\frac{-3}{1}]^{i}[\frac{-3}{2}]$ for all $0\leq i\leq 1_{min,A}-1$, \\
$\Join^{(0)}:=\bigstar_{1_{min,A}-2}$,\\
$\alpha_i:=(-1)^{(2_{min,A}+1)}\ ^{(2_{min,A}-1)}C_i$ for each $0\leq i\leq 2_{min,A}-1$,\\ 
$\beta_j^{(0)}:=\alpha_{j-1}+\alpha_{j}$ for each $1\leq j\leq 1_{min,A}-1$ and $\beta_0^{(0)}:=\alpha_0$.\\
\end{center}
In each column, the subscript $j$ of $\beta_j^{(0)}$ decreases (by $1$ in each consecutive row) as one goes down.\\

The system $\mathfrak{B}_{A,f}^{eq}\mathfrak{Y}^{eq}=\mathfrak{e}^{eq}$ of linear equations is solvable if and only if the rank of the matrix $\mathcal{C}_{A,f}$ equals the rank of the matrix $[\mathcal{C}_{A,f}\ \mathcal{E}]$ where $\mathcal{E}$ is the vector $(0,\ldots,0,1)^t$. But unlike in subcase (II.1), it is not clear here whether or not the ranks of these two matrices are equal. Hence one needs to compute the ranks of such matrices to determine whether or not the system $\mathfrak{B}_{A,f}^{eq}\mathfrak{Y}^{eq}=\mathfrak{e}^{eq}$ is solvable.\\

\noindent \underline{Subcase (II.3)}: When $(-3)_{max,A}<max\{1_{min,A},2_{min,A}\}$.\\

\noindent Possibility (a): When $(-3)_{max,A}<2_{min,A}$ but $(-3)_{max,A}\geq 1_{min,A}$. \\
The matrix $\mathfrak{B}_{A,f}^{tr}$ consists only of Part A, there is no Part B. The structure of Part A of the matrix $\mathfrak{B}_{A,f}^{tr}$ is similar to that in subcase(II.2). Hence the matrix $\mathfrak{B}_{A,f}^{eq}$ consists only of Part A, there is no Part B. The system $\mathfrak{B}_{A,f}^{eq}\mathfrak{Y}^{eq}=\mathfrak{e}^{eq}$ of linear equations is not solvable because $1_{min,A}\geq 1$ for any monomial $A$ and therefore the last row of the matrix corresponding to Part A of $\mathfrak{B}_{A,f}^{eq}$ is the zero row.\\

\noindent Possibility (b): When $(-3)_{max,A}<2_{min,A}$ and $1_{min,A}$ both. \\
Neither Part A nor Part B of the matrix $\mathfrak{B}_{A,f}^{tr}$ exists. So the system $\mathfrak{B}_{A,f}^{tr}\mathfrak{Y}^{tr}=\mathfrak{e}^{tr}$ of linear equations is not solvable by default. \\

\noindent Possibility (c): When $(-3)_{max,A}\geq 2_{min,A}$ and $(-3)_{max,A}<1_{min,A}$. \\
The matrix $\mathfrak{B}_{A,f}^{tr}$ consists only of Part B, there is no Part A. Moreover, the last row of the matrix $\mathfrak{B}_{A,f}^{tr}$ (which consists only of Part B) is the zero row. Hence the system $\mathfrak{B}_{A,f}^{tr}\mathfrak{Y}^{tr}=\mathfrak{e}^{tr}$ of linear equations is not solvable because the last component of the vector $\mathfrak{e}^{tr}$ is non-zero.\\
\section{The algorithm and suspicion about irrationality}\label{s.the-algorithm-irrational}
Recall from section \ref{s.statetheproblem} that we are interested in computing the Hilbert-Kunz function of $R$ where $R$ is a `disjoint-term trinomial hypersurface'. Over any field of characteristic $p>0$, the value of this Hilbert-kunz function at $p^n$ is given by the length $l(\frac{R}{{\mathfrak{m}}^{(p^n)}+J})$ for any given positive integer $n$. It follows from subsection \ref{ss.filtration-and-checking} that this length equals the cardinality of the set $\{A\in\mathfrak{M}|A\notin A_c+J\}$ where the set $\mathfrak{M}$ is as defined in \S 2.1 of \cite{Su} and the ideals $A_c$ and $J$ are as mentioned in subsection \ref{ss.mutation} and section \ref{s.statetheproblem} respectively. It is therefore equivalent to compute the cardinality of the set $\{A\in\mathfrak{M}|A\in A_c+J\}$. Recall that $f=[3]+[2]+[1]$. Given any monomial $A\in\mathfrak{M}$, the following is an algorithm to check whether or not $A\in A_c+J$:\\
\begin{itemize}
 \item Compute $1_{min,A},2_{min,A},(-2)_{max,A}$ and $(-3)_{max,A}$ for the monomial $A$ (see \S \ref{n.partition-of-L(A,f)} for the meaning of the notation $1_{min,A},2_{min,A},(-2)_{max,A}$ and $(-3)_{max,A}$.).
 \item If condition (i) of theorem \ref{t.mutation-stop-modifiedv4} holds, declare that $A\in A_c+J$.
 \item If condition (ii) of theorem \ref{t.mutation-stop-modifiedv4} holds, declare that $A\in A_c+J$.
 \item If neither condition (i) nor (ii) of theorem \ref{t.mutation-stop-modifiedv4} holds but $[3]$ divides $A$, then check that $(-2)_{max,A}=0$ or not.
\end{itemize}
\noindent If $(-2)_{max,A}=0$, then proceed in the following way:\\
\begin{itemize}
  \item If $(-3)_{max,A}\geq 1_{min,A}+2_{min,A}-1$, declare that $A\in A_c+J$.
  \item If $max\{1_{min,A},2_{min,A}\}\leq (-3)_{max,A}<1_{min,A}+2_{min,A}-1$, it is not clear whether $A\in A_c+J$ or not. For knowing whether $A\in A_c+J$ or not, we need to do a rank computation corresponding to matrices of the type given in table $8$ of subcase (II.2) of section \ref{s.solvability-of-thenewsys}.  
  \item If $(-3)_{max,A}<max\{1_{min,A},2_{min,A}\}$, declare that $A\notin A_c+J$.
\end{itemize}
\noindent If $(-2)_{max,A}\geq 1$, then proceed in the following way:\\
\begin{itemize}
 \item If $(-3)_{max,A}\geq 1_{min,A}+2_{min,A}-1$, declare that $A\in A_c+J$.
 \item If $max\{1_{min,A},2_{min,A}\}\leq (-3)_{max,A}<1_{min,A}+2_{min,A}-1$, it is not clear whether $A\in A_c+J$ or not. For knowing whether $A\in A_c+J$ or not, we need to do a rank computation corresponding to matrices of the type given in table $5$ of subcase (I.2) of section \ref{s.solvability-of-thenewsys}.  
 \item If $(-3)_{max,A}<2_{min,A}$ but $(-3)_{max,A}\geq 1_{min,A}$, then $A\in A_c+J$ if and only if $1_{min,A}\leq (-2)_{max,A}$.
 \item If $(-3)_{max,A}<2_{min,A}$ and $(-3)_{max,A}<1_{min,A}\leq (-2)_{max,A}+(-3)_{max,A}$, then $A\in A_c+J$ if and only if $1_{min,A}\leq (-2)_{max,A}$.
 \item If $(-3)_{max,A}<2_{min,A}$ and $(-3)_{max,A}<1_{min,A}$ but $1_{min,A}>(-2)_{max,A}+(-3)_{max,A}$, declare that $A\notin A_c+J$. 
 \item If $(-3)_{max,A}\geq 2_{min,A}$ and $(-3)_{max,A}<max\{1_{min,A},2_{min,A}\}$, it is not clear whether $A\in A_c+J$ or not. For knowing whether $A\in A_c+J$ or not, we need to do a rank computation corresponding to matrices of the type given in table $6$ of subcase (I.3) of section \ref{s.solvability-of-thenewsys}.  
\end{itemize}
In this algorithm, there are certain situations where it is not clear whether $A\in A_c+J$ or not. In such situations, one needs to compute ranks of some huge matrices of the types mentioned above and the existence of such weird situations make me suspect that the Hilbert-Kunz multiplicity for some `disjoint-term trinomial hypersurfaces' can become irrational. More work regarding such rank computation is in progress. A high end computing platform may help solve particular examples to a great extent, but it cannot give a formula for the Hilbert-Kunz function. The matrices whose ranks need to be computed are of the types mentioned in tables $5,6$ and $8$ above. These matrices exhibit a nice combinatorial pattern, which in turn can show us a rhythm in which ranks of the concerned matrices are going to vary as we move from one monomial to another in the set $\mathfrak{M}$. This rhythm will help us provide a concrete reasoning to my suspicion regarding the irrationality of certain Hilbert-Kunz multiplicities as mentioned above. I wish to address this problem at least over fields of characteristic $2$ in my next work. 
\providecommand{\bysame}{\leavevmode\hbox
to3em{\hrulefill}\thinspace}
\providecommand{\MR}{\relax\ifhmode\unskip\space\fi MR }
\providecommand{\MRhref}[2]{%
  \href{http://www.ams.org/mathscinet-getitem?mr=#1}{#2}
} \providecommand{\href}[2]{#2}


\begin{thebibliography}{10}

\bibitem{Mo}
Paul Monsky, \emph{The Hilbert-Kunz function}, math, Ann. 263 (1983), 43--49.

\bibitem{Mo4}
Paul Monsky, \emph{The Hilbert-Kunz multiplicity of an irreducible trinomial}, Journal of Algebra, 304 (2006), 1101--1107

\bibitem{Buch}
Ragnar-Olaf Buchweitz, Qun Chen, \emph{Hilbert-Kunz functions of Cubic Curves and Surfaces}

\bibitem{Conca}
A. Conca, \emph{Hilbert-Kunz function of monomial ideals and binomial hypersurfaces}, Manuscripta Math. 90 (1996), 287--300

\bibitem{Trivedi}
V. Trivedi, \emph{Semistability and Hilbert-Kunz multiplicities for curves}, Journal of Algebra, 284 (2005), 627--644

\bibitem{Mo1}
Paul Monsky, \emph{Rationality of Hilbert-Kunz multiplicities: A likely Counterexample}, Michigan Math.J.57 (2008)

\bibitem{Mo2}
Paul Monsky, \emph{Transcendence of some Hilbert-Kunz multiplicities (modulo a conjecture)}, arXiv:0908.0971v1, [math.AC]

\bibitem{Mo3}
Paul Monsky, \emph{Algebraicity of some Hilbert-Kunz multiplicities (modulo a conjecture)}, arXiv:0907.2470v1, [math.AC]

\bibitem{Su}
Shyamashree Upadhyay, \emph{The Hilbert-Kunz function for Binomial Hypersurfaces}, arXiv:1101.5936, [math.CO]

\bibitem{Han}
C. Han, \emph{The Hilbert-Kunz function of a diagonal hypersurface}, PhD thesis, Brandeis University, 1991, MR2635348.
\end{thebibliography}
\end{document}